%% file: Markov_Diagram.tex
\documentclass[10pt,a4paper,twoside]{amsart}

\usepackage[utf8]{inputenc}
\usepackage[english]{babel}
\usepackage[T1]{fontenc}
\usepackage{microtype, fancyhdr, lmodern, xcolor}

\usepackage[textsize=9]{todonotes}

\usepackage{amsfonts, amsmath, amsthm, amssymb, mathrsfs, amscd}		
\numberwithin{equation}{section}	
\usepackage{faktor}
\allowdisplaybreaks

\usepackage{url, enumitem}
\usepackage[noadjust]{cite}

\usepackage{multirow,bigdelim, caption}
\usepackage{booktabs}

\usepackage{graphicx}

\usepackage{hyperref}

\theoremstyle{plain}
\newtheorem{thm}{Theorem}[section]
\newtheorem{prop}[thm]{Proposition}
\newtheorem{cor}[thm]{Corollary}
\newtheorem{lem}[thm]{Lemma}

\theoremstyle{remark}
\newtheorem{rmk}[thm]{Remark}

\theoremstyle{definition}
\newtheorem{defn}[thm]{Definition}
\newtheorem{exam}[thm]{Example}

\input{Macro}

\begin{document}

\title{A Markov's theorem for extended welded braids and links}

\author[Damiani]{Celeste Damiani}
\address{Department of Mathematics,
Osaka City University,
Sugimoto, Sumiyoshi-ku,
Osaka 558-8585, Japan}
\email{celeste.damiani@math.cnrs.fr}

\subjclass[2010]{Primary 20F36; secondary 57Q45}

\keywords{Braid groups, welded braids, ribbon braids, welded links, ribbon torus-links}

\date{\today}

\begin{abstract}
Extended welded links are a generalization of Fenn, Rim\'{a}nyi, and Rourke's welded links. Their braided counterpart are extended welded braids, which are closely related to ribbon braids and loop braids. In this paper we prove versions of Alexander and Markov's theorems for extended welded braids and links, following Kamada's approach to the case of welded objects.
\end{abstract}

\maketitle

\section{Introduction}

%
%
\subsection*{State of the art}
\emph{Welded links} were introduced by Fenn-Rim\'{a}nyi-Rourke~\cite{Fenn-Rimanyi-Rourke:1997} as equivalence classes of link diagrams in the $2$-dimensional space. They can be considered as virtual links up to additional Reidemeister moves called \emph{forbidden moves}. Satoh \cite{Satoh:2000} considered the relation between welded links and \emph{ribbon torus-links}. He extended a construction of Yajima~\cite{Yajima:2manifolds}, defining a surjective map $Tube$ from welded knots to ribbon torus-knots, which allows to associate to any ribbon torus-knot a welded knot. 
This fact suggested that ribbon torus-knots could be the topological counterparts of welded knots. However, the $Tube$ map is not injective: for instance, it is invariant under the \emph{horizontal mirror image} on welded diagrams \cite[Proposition~3.3]{Ichimori-Kanenobu:2012} (see also \cite{Winter:2009, Satoh:2000}), while welded links in general are not equivalent to their horizontal mirror images. 

Their braided counterparts are \emph{welded braids}, also introduced in~\cite{Fenn-Rimanyi-Rourke:1997}. Welded braid groups can be seen as quotients of virtual braid groups. Fenn-Rim\'{a}nyi-Rourke proved that these groups are isomorphic to the groups of \emph{braid-permutation automorphisms} of the free groups. They are also isomorphic to groups appearing in many other contexts, see for instance the surveys~\cite{Damiani:Journey, BarNatan-Dancso:Survey}. In particular, they are isomorphic to the groups of \emph{ribbon braids}, which are the braided counterparts of ribbon torus-links. In fact, on braided objects, the $Tube$ map is an isomorphism~\cite{Audoux-Bellingeri-Al:HomotopyRibbonTubes}. 

For welded braids and links we have versions of Alexander's and Markov's theorems, due to Kauffman-Lambropoulou~\cite{Kauffman-Lambropoulou:L-move} and to Kamada~\cite{Kamada:Markov}. The isomorphism between welded braid groups and ribbon braid groups guarantees that Alexander's theorem for welded objects holds when passing to ribbon braids and ribbon torus-links. However, the lack of a bijection between welded links and ribbon torus-links impedes us to translate Markov's theorem for welded objects to ribbon braids and ribbon torus-links.

\subsection*{Motivation and contribution}
The aim of this note is to make a step towards a Markov's theorem for ribbon braids and ribbon torus-links. This is done by studying a class of objects that appear as suitable candidates to be a diagrammatical representation of ribbon torus-links, while remaining in the domain of usual link diagrams. In fact ribbon torus-links can be also represented by chord-diagrams, see for instance~\cite{Kawauchi:Chord16moves}. The objects we consider are an enhanced version of welded links: they are called \emph{extended welded links} and have been  introduced in~\cite{Damiani:Journey}. In a certain sense, they can be seen as a quotient of welded links: in fact each extended welded link is equivalent to a welded link, but the introduction of certain marks on the diagrams, called \emph{wen marks}, adds some moves to the generalized Reidemeister relations, making it possible that two non-equivalent welded links, when considered as extended welded links, become equivalent through moves involving wen marks. 

There are two reasons that point to extended welded links as promising candidates to being diagrammatical representations for ribbon torus-links. The first reason is that that extended welded link diagrams are equivalent to their sign reverse, and in consequence to their  horizontal mirror image~(Proposition~\ref{P:Mirror}). As for the second reason, let us consider \emph{extended ribbon braids} by allowing \emph{wens}, which are embeddings in the $4$-dimensional space of a Klein bottle cut along a meridional circle, on the braided annuli that compose ribbon braids. The groups of extended ribbon braids appear when looking for a version of Markov theorem for ribbon braids and torus-links in $B^3 \times S^1$. In fact it can be proven that taken a pair of ribbon braids, their closures are isotopic as ribbon torus-links in $B^3 \times S^1$ if and only if they are conjugate as \emph{extended} ribbon braids~\cite{Damiani:Tesi}. However, if one considers extended ribbon braids to begin with, the statement is the exact analogue of the usual case in dimension~$3$: taken a pair of extended ribbon braids, their closures are isotopic as ribbon torus-links in $B^3 \times S^1$ if and only if they are conjugate as extended ribbon braids. 
This is relevant to this paper because the groups of extended ribbon braids  are isomorphic to the groups of \emph{extended welded braids}, which are the braided counterpart of extended welded links~\cite[Theorem 6.12]{Damiani:Journey}. 

The main result of this paper is the following:
\newtheorem*{thm:MarkovIntro}{Theorem \ref{T:Markov}}
\begin{thm:MarkovIntro}
Two ext. welded braid diagrams that admit closure have equivalent closures as ext. welded link diagrams if and only if they are related by a finite sequence on the following moves:
\begin{description}[before={\renewcommand\makelabel[1]{\bfseries ##1}}]
\item[$(M0)$] isotopy of $\Rr^2$ and generalized  Reidemeister moves;
\item[$(M1)$] conjugation in the welded ext. braid group;
\item[$(M2)$] a right stabilization of positive, negative or welded type, and its inverse operation.
\end{description}
\end{thm:MarkovIntro}

%

%

%
\subsection*{Structure of the paper}
In Section~\ref{S:ExtWeldBraids} we introduce extended welded braid diagrams and the extended welded braid groups. We give a presentation for them and describe their relation with virtual and welded braids. In Section~\ref{S:ExtWeldLinks} we discuss extended welded links and give a combinatorial description for them in terms of Gauss data. We state a version of Alexander's Theorem for extended welded objects (Proposition~\ref{P:Alexander}) and state some results that allow us to use Gauss data to describe extended welded links.
In Section~\ref{S:MarkovTheorem} we prove the main result (Theorem~\ref{T:Markov}). Finally, in Section~\ref{S:Mirror} we show that extended welded knots are equivalent to their horizontal mirror images (Proposition~\ref{P:Mirror}).

\section{Extended welded braids}
\label{S:ExtWeldBraids}
\begin{figure}[hbt]
\centering
\includegraphics[scale=0.6]{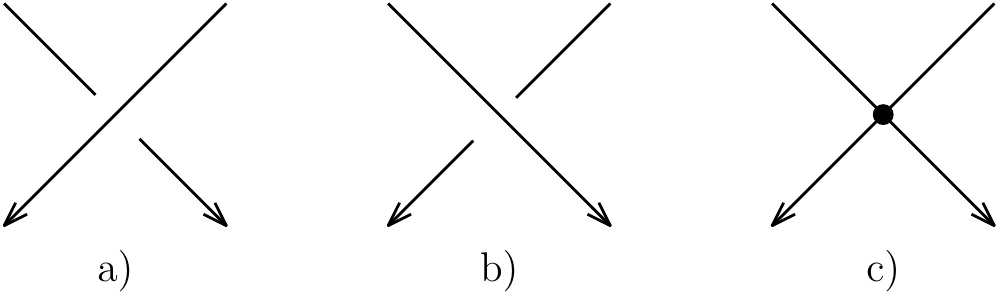}
\caption{a) Classical positive crossing, b) Classical negative crossing, c) Welded crossing.}
\label{F:Crossings}
\end{figure}

\begin{figure}[htb]
\centering
\includegraphics[scale=0.6]{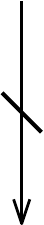}
\caption{A wen mark on a strand.}
\label{F:WenMark}
\end{figure}

An \emph{extended welded braid diagram}, or \emph{ext. welded braid diagram} on $\nn$ strings is a planar diagram composed by a set of $\nn$ oriented and monotone $1$-manifolds immersed in $\Rr^2$ starting from $\nn$ points on a horizontal line at the top of the diagram down to a similar set of $\nn$ points at the bottom of the diagram. The $1$-manifolds are allowed to cross in transverse double points, which will be decorated in three kinds of ways, as shown in Figure~\ref{F:Crossings}. Depending on the decoration, double points will be called: \emph{classical positive} crossings, \emph{classical negative} crossings and \emph{welded} crossings. On each $1$-manifold there can possibly be marks as in~Figure~\ref{F:WenMark}, which we will call \emph{wen marks}.

\begin{figure}[hbt]
\centering
\includegraphics[scale=0.6]{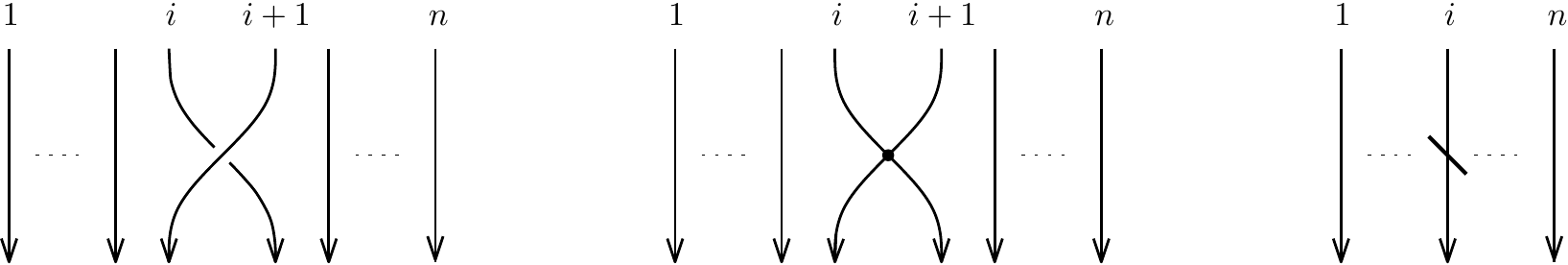}
\caption{Elementary diagrams $\sig\ii$, $\rr\ii$, and~$\tau_\ii$.}
\label{F:SigRhoTau}
\end{figure}

Let us assume that the double points occur at different $y$-coordinates. Then an ext. welded braid diagram determines a word in the elementary diagrams illustrated in Figure~\ref{F:SigRhoTau}. We call $\sig\ii$ the elementary diagram representing the $(i+1)$-th strand passing over the $i$-th strand, $\rr\ii$ the welded crossing of the strands $i$ and~$(i+1)$, and $\tau_\ii$ the wen mark diagram.

\begin{figure}[hbt]
\centering
\includegraphics[scale=.6]{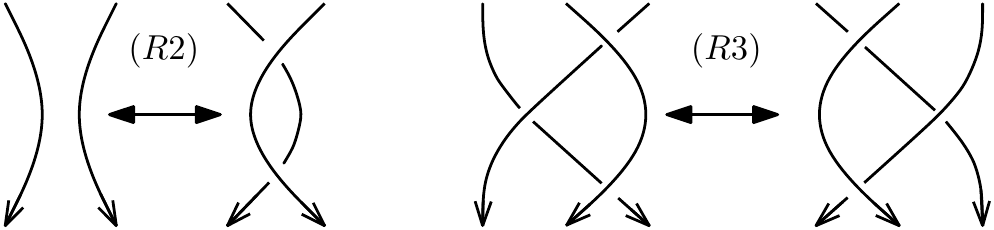}
\caption{Classical Reidemeister moves for braid-like objects.}
\label{F:Classical}
\end{figure}

\begin{figure}[hbt]
\centering
\includegraphics[scale=.6]{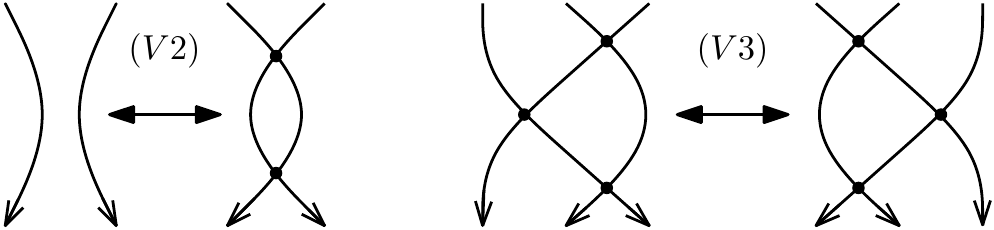}
\caption{Virtual Reidemeister moves for braid-like objects.}
\label{F:Virtual}
\end{figure}

\begin{figure}[hbt]
\centering
\includegraphics[scale=.6]{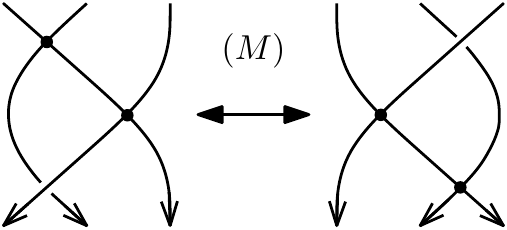}
\caption{Mixed Reidemeister moves.}
\label{F:Mixed}
\end{figure}

\begin{figure}[hbt]
\centering
\includegraphics[scale=.6]{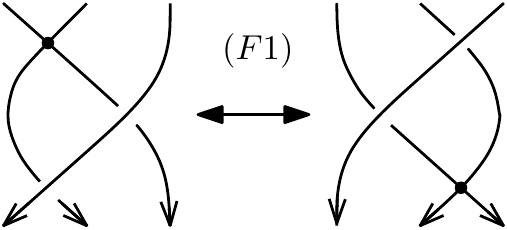}
\caption{Welded Reidemeister moves.}
\label{F:Welded}
\end{figure}

\begin{figure}[htb]
\centering
\includegraphics[scale=0.6]{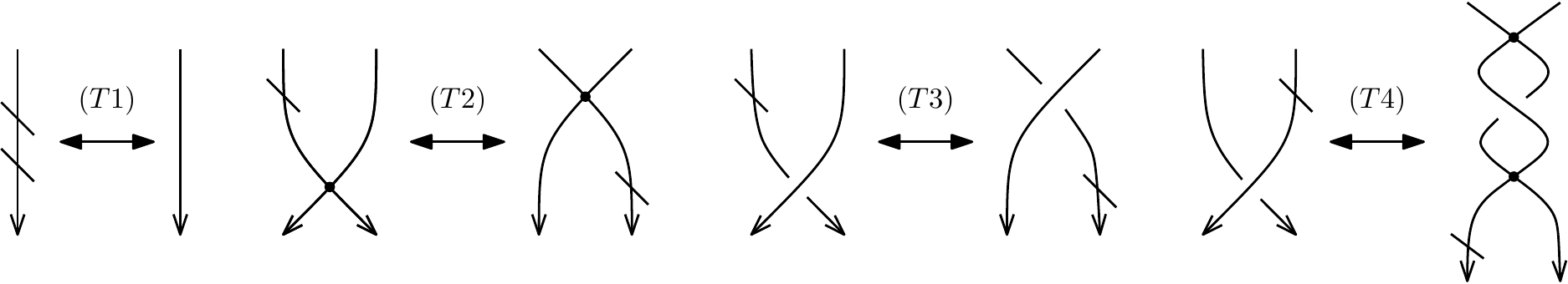}
\caption{Extended Reidemeister moves.}
\label{F:ReidWen}
\end{figure}

\begin{defn}
\label{D:bWelded}
An \emph{extended welded braid}, or \emph{ext. welded braid} to keep notation short, is an equivalence class of ext. welded braid diagrams under the equivalence relation given by isotopy of $\Rr^2$ and the following moves: 
\begin{itemize}
\item classical Reidemester moves (Figure~\ref{F:Classical});
\item virtual Reidemeister moves (Figure~\ref{F:Virtual});
\item mixed Reidemeister moves (Figure~\ref{F:Mixed});
\item welded Reidemeister moves (Figure~\ref{F:Welded});
\item extended Reidemester moves (Figure~\ref{F:ReidWen}).
\end{itemize} 
This equivalence relation is called \emph{(braid) generalized Reidemeister equivalence}. 
For~$\nn \geq 1$, the \emph{extended welded braid group} on $\nn$ strands $\WBE\nn$ is the group of equivalence classes of ext. welded braid diagrams by generalized Reidemeister equivalence. The group structure on these objects is given by: stacking and rescaling as product, braid mirror image as inverse, and the trivial diagram as identity.
\end{defn}

\begin{rmk}
If wen marks were not allowed, the group defined would be the group of \emph{welded braids} $\WB\nn$, introduced by Fenn, Rim\'{a}nyi and Rourke in~\cite{Fenn-Rimanyi-Rourke:1997}. This group is isomorphic to \emph{loop braid groups} $\LB\nn$, \emph{ribbon braid groups} $\rB\nn$, and many others. More precisely, we can see welded braids as the diagrams describing loop braids. For a survey on these groups, see~\cite{Damiani:Journey}. For more details about loop braid groups seen in different contexts, this is a non-exhaustive list of references:~ \cite{BaezWiseCrans:Exotic, BrendleHatcher:2013, Goldsmith:MotionGroups, Lin:2008,   McCool:1986, Savushkina:1996, Wattenberg:1972 }.
\end{rmk}

Ext. welded braid groups $\WBE\nn$ are isomorphic to \emph{ring groups $R_n$} discussed in \cite{BrendleHatcher:2013} (see \cite{Damiani:Journey} for a proof of the equivalence). Therefore, they admit a presentation given by generators $\{\sig\ii, \rr\ii \mid \ii=1, \ldots, \nno\}$ and $\{\tau_\ii \mid \ii=1, \dots , \nn \}$, subjected to the following relations:
\begin{equation}
\label{E:Rpresentation}
\begin{cases}
\sig{i} \sig j = \sig j \sig{i}  \, &\text{for } \vert  i-j\vert > 1\\
\sig{i} \sig {i+1} \sig{i} = \sig{i+1} \sig{i} \sig{i+1} \, &\text{for } i=1, \dots, \nn-2 \\
\rr{i} \rr j = \rr j \rr{i}  \, &\text{for }  \vert  i-j\vert > 1\\
\rr\ii\rr{i+1}\rr\ii = \rr{i+1}\rr\ii\rr{i+1} \, &\text{for }  i=1, \dots, \nn-2  \\
\rrq{i}2 =1 \, &\text{for }  i=1, \dots, \nno \\
\rr{i} \sig{j} = \sig{j} \rr{i}   \, &\text{for }  \vert  i-j\vert > 1\\
\rr{i+1} \rr{i} \sig{i+1} = \sig{i} \rr{i+1} \rr{i} \,  &\text{for }  i=1, \dots, \nn-2  \\
\sig{i+1} \sig{i} \rr{i+1} = \rr{i} \sig{i+1} \sig{i}  \,  &\text{for }  i=1, \dots, \nn-2, \\
\tau_{i} \tau_j = \tau_j \tau_{i}  \, &\text{for }    i \neq j \\
\tau_\ii^2=1 \, &\text{for }  i=1, \dots, \nn \\
\sig{i} \tau_j = \tau_j \sig{i}  \, &\text{for }  \vert  i-j\vert > 1 \\
\rr{i} \tau_j = \tau_j \rr{i}  \, &\text{for }  \vert  i-j\vert > 1 \\
\tau_i \rr\ii = \rr\ii \tau_{i+1} \, &\text{for }  i=1, \dots, \nno  \\
\tau_\ii \sig\ii = \sig\ii \tau_{i+1}  \, &\text{for } i=1, \dots, \nno \\
\tau_{i+1} \sig\ii = \rr\ii \siginv\ii \rr\ii \tau_\ii  \, &\text{for } i=1, \dots, \nno. \\
\end{cases}
\end{equation}

\subsection{Extended Markov moves}
\label{SS:Extended_Markov_Moves}
\begin{figure}[htb]
\centering
\includegraphics[scale=0.6]{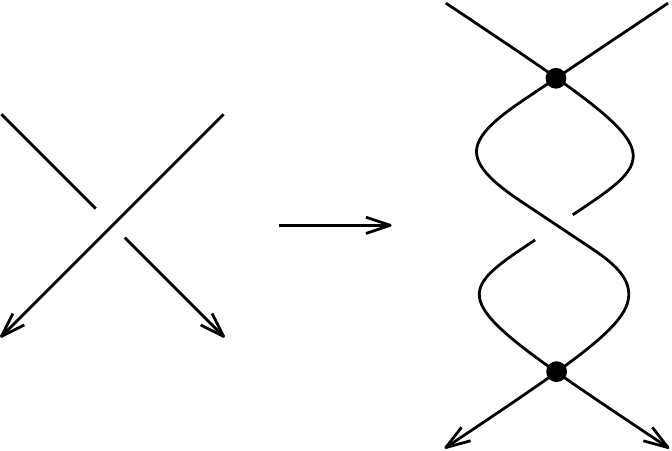}
\caption{Sign reversal.}
\label{F:Sign_Reversal}
\end{figure}

We introduce here some ``moves'' on ext. welded braid diagrams. Let $b_1$ and $b_2$ be two ext. welded braid diagrams. 
\begin{itemize}
\item If $b_1$ and $b_2$ are represent the same ext. welded braid, then we say that $b_2$ if obtained from $b_1$ by a  \emph{$(M0)$-move}.
\item If $b_1$ and $b_2$ have the same degree (the same number of strands), then the ext. welded braid diagram $b_1 b_2$ is obtained from the diagram $b_2 b_1$ by \emph{conjugation}, also called \emph{$(M1)$-move}.
\item Suppose $b_1$ has degree $\nn$. Then $\iota(b_1)$ is the diagram of degree $\nn +1$ obtained by adding a trivial strand to the right of $b_1$. If $b_2 = \iota(b_1)\sig\nn$, or $\iota(b_1)\siginv\nn$, or $\iota(b_1)\rr\nn$, we say that $b_2$ is obtained from $b_1$ by \emph{ positive, negative, or welded stabilization}, also called \emph{$(M2)$-move}.
\end{itemize}
We introduce also a third type of move, which is a direct consequence of moves $(M0)$ and $(M1)$, as we will prove in Proposition~\ref{P:sign_reversal}, but will be used in Section~\ref{S:Mirror} to enlighten an important property that extended welded links have and (non extended) welded links don't. 
\begin{itemize}
\item If $b_2$ is obtained from $b_1$ by replacing each classical crossing of $b_1$ with a crossing of the opposite sign conjugated by welded crossings as in Figure~\ref{F:Sign_Reversal}, we say that $b_2$ is obtained by \emph{sign reversal}. 
\end{itemize}

\section{Extended welded links and Gauss data}
\label{S:ExtWeldLinks}
\subsection{Extended welded links}

An \emph{extended welded link diagram}, or just \emph{ext. welded link diagram}, is the immersion in~$\Rr^2$ of a collection of disjoint, closed, oriented $1$-manifolds such that all multiple points are transverse double points. Double points are decorated with classical positive, classical negative, or welded information as in Figure~\ref{F:Crossings}. On each $1$-manifod there can possibly be an even number of wen marks as in~Figure~\ref{F:WenMark}. We assume that ext. welded link diagrams are the same if they are isotopic in~$\Rr^2$. 
Taken an ext. welded link diagram $K$, we call \emph{real crossings} its set of classical positive and classical negative crossings.

\begin{figure}[htb]
	\centering
		\includegraphics[scale=0.7]{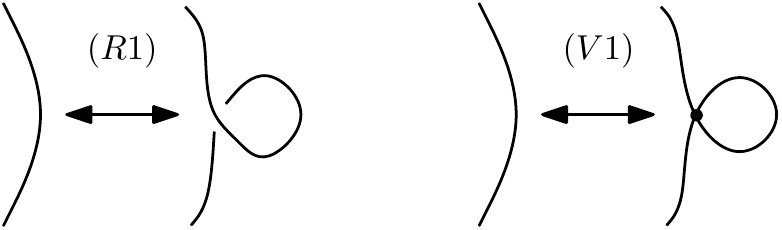} 
	\caption{Reidemeister moves of type I.}
	\label{F:ReidOne}
\end{figure}

\begin{defn}
\label{D:ExWeldedLink}
An \emph{extended welded link} or \emph{ext. welded link} is an equivalence class of ext. welded link diagrams under the equivalence relation given by isotopies of $\Rr^2$, moves from Definition~\ref{D:bWelded}, and by classical and virtual Reidemeister moves $(R1)$ and $(V1)$ as in Figure~\ref{F:ReidOne}. This equivalence relation is called \emph{generalized Reidemeister equivalence}. 
\end{defn}

\begin{rmk}
The reason because wen marks on ext. welded links can only appear with even parity on each component of a ext. welded link diagram is motivated by the relation between ext. welded objects and ext. ribbon objects. For more details, see~\cite[proof of Proposition 2.4]{Audoux:2016}.
\end{rmk}

\begin{figure}[htb]
\centering
\includegraphics[scale=0.6]{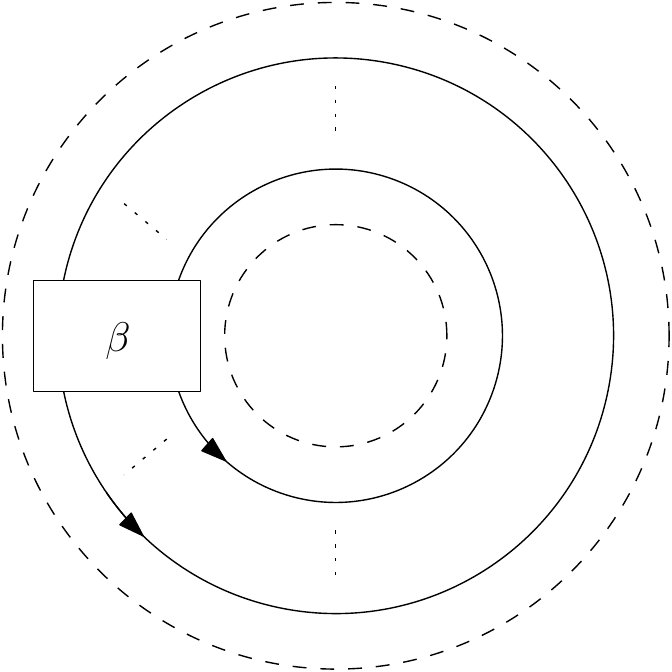}
\caption{Closure of an ext. welded braid diagram.}
\label{F:closure}
\end{figure}

The \emph{closure} of an ext. welded braid diagram is obtained as for usual braid diagrams (see Figure~\ref{F:closure}), with the condition that ext. welded braids can be closed only when they have an even number of wen marks on each component. 

\begin{figure}[htb]
\centering
\includegraphics[scale=0.61]{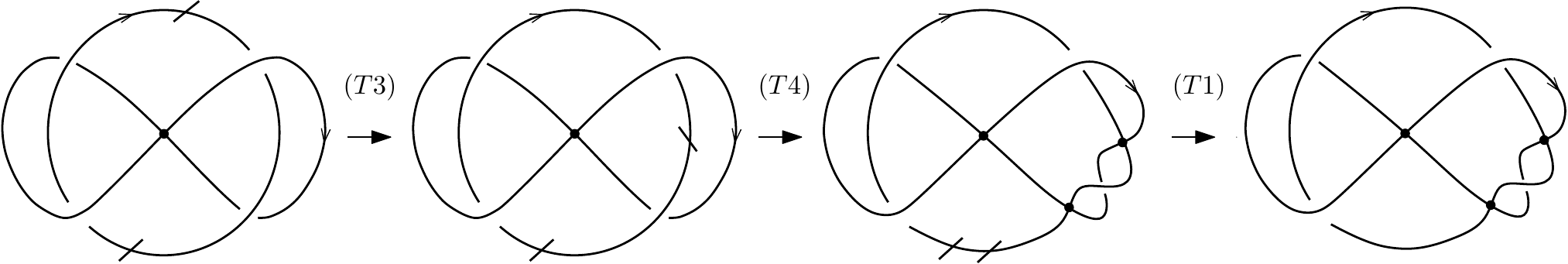}
\caption{Sliding wen marks along a diagram in order to obtain a wen marks-free diagram.}
\label{F:wen_slide}
\end{figure}

\begin{prop}
\label{P:Alexander}
Any ext. welded link can be described as the closure of an ext. welded braid diagram which is generalized Reidemeister equivalent to a (non-extended) welded braid diagram. 
\end{prop}
\proof
Recall that an ext. welded link can have only an even number of wen marks on each component. Then, taken $l$ to be an ext. welded link, and $L$ an ext. welded link diagram representing it, it is always possible to find a diagram $L'$ without wen marks. This is done by making one wen mark slide along the component it belongs to, in order to make it adjacent to another wen mark, and the cancelling the wen marks pairwise. For an example, see Figure~\ref{F:wen_slide}. 
Then it is enough to prove that any welded link can be described as the closure of a welded braid. This is done in \cite[Proposition~8]{Kamada:Markov} and \cite[Theorem~1]{Kauffman-Lambropoulou:L-move}.

Alternatively, this result can be proved directly by giving a braiding algorithm, as the one we present in Subsection~\ref{SS:Gauss_Data}.
\endproof

\subsection{Gauss data}
\label{SS:Gauss_Data}

\begin{figure}[htb]
\centering
\includegraphics[scale=0.6]{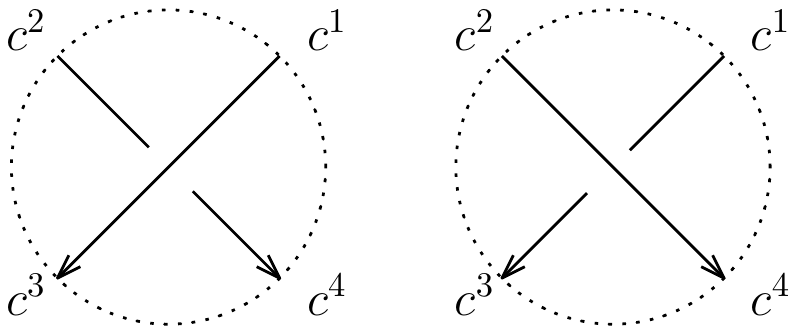}
\caption{Intersection points $c^1, c^2, c^3$ and $c^4$ composing $\partial N(c) \cap K$.}
\label{F:Crossing_Neighbourhood}
\end{figure}

We can associate to an ext. welded knot diagram a \emph{Gauss data}. This is a combinatoric description of a Gauss diagram, as intended in the spirit of~\cite{Goussarov-Polyak-Viro:2000}. In the context of ext.welded diagrams, Gauss data needs to contain the added information for wen marks. We recall and adapt here the description of Gauss data given in~\cite{Kamada:Markov}. Let $K$ be an ext. welded link diagram. We introduce some  notation.
\begin{itemize}
\item We denote by $C_K$ the set of positive and negative crossings of $K$, also called the set of \emph{real crossings}.
\item We define a sign map  $S_K \colon C_K \to \{ -1, 1\}$  on the set real crossings, sending positive crossings to $1$ and negative crossings to~$-1$.
\item For a real crossing $c \in C_K$, let $N_c$ be a regular neighbourhood of $c$; we denote by $W_K$ the closure of $\Rr^2 \setminus \bigcup_{c \in C_K} N_c$,  and by $K\restriction_{W_K}$ the restriction of $K$ to~$W_K$.
\item We denote by $c^1, c^2, c^3$ and $c^4$ the four points that compose $\partial N(c) \cap K$, and by $C_K^\partial$ the set $\{ c^\ii \mid c \in C_K,  \ii\in \{1, \ldots, 4 \} \}$. See Figure~\ref{F:Crossing_Neighbourhood}. 
\end{itemize}

 Define a subset $G_K \subset C_K^\partial \times C_K^\partial \times \faktor{\Zz}{2\Zz}$ such that $(a, b, n+2\Zz) \in G_K$ if and only if the restriction of $K$ to $W_K$ has an oriented arc starting from $a$ and terminating at~$b$, and $\nn$ is the number of wen marks on it. Elements $(a, b, 0+2\Zz)$ and $(a, b, 1+2\Zz)$ will be respectively denoted by $(a, b)$ and~$\overline{(a, b)}$. Finally, we denote by $\mu_K$ the number of components of~$K$.

\begin{defn}
\label{D:Gauss_data}
The Gauss data of $K$ is the quadruple $G(K) = ( C_K, S_K, G_K, \mu_K)$. We say that two Gauss ext. welded link diagrams $K$ and $K^\prime$ have the \emph{same Gauss data} if the have the same number of components $\mu_K=\mu_{K^\prime}$ and if there is a bijection $g \colon C_K \to C_{K^\prime}$ such that $g$ preserves the signs of the crossings, and the presence of welded marks. This means that if $(a, b)$ is in $G_K$, then $(g(a), g(b))$ is in  $G_{K^\prime}$, and  if $\overline{(a, b)}$ is in $G_K$, then $\overline{(g(a), g(b))}$ is in~$G_{K^\prime}$.
\end{defn}

\begin{figure}[htb]
\centering
\includegraphics[scale=0.8]{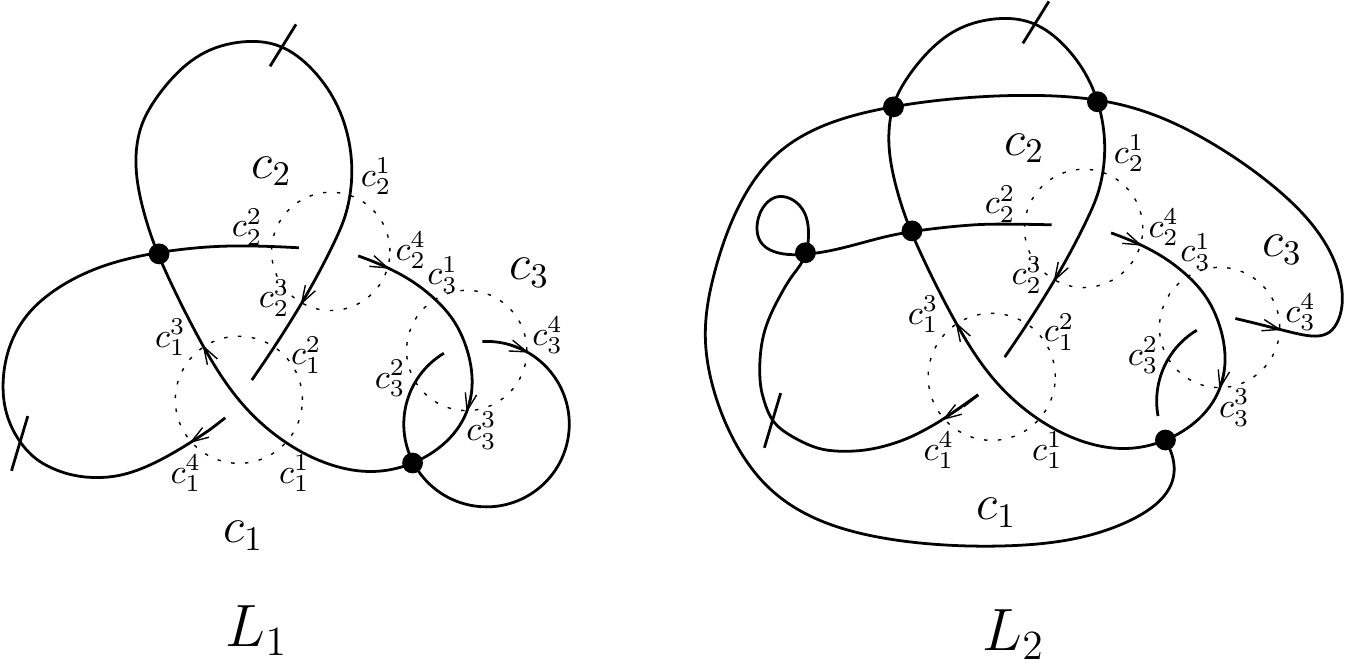}
\caption{Two ext. welded link diagram with the same Gauss data.}
\label{F:Trifoglio_Gauss}
\end{figure}

\begin{exam}
Let us compute the Gauss data for link $L_1$ in Figure~\ref{F:Trifoglio_Gauss}. Its Gauss data is given by:
\begin{itemize}
\item $C_{L_1} = \{c_1, c_2, c_3\}$;
\item $S_{L_1} = c_1 \longmapsto +1, c_2 \longmapsto +1, c_3 \longmapsto +1$;
\item $G_{L_1} = \big\{ (c_2^3, c_1^2), \overline{(c_1^4, c_2^2)}, (c_2^4, c_3^1), (c_3^3, c_1^1), \overline{(c_1^3, c_2^1)}, (c_3^4, c_3^2) \big\}$;
\item $\mu_{L_1} = 2$.
\end{itemize}

\end{exam}

Let $K$ be an ext. welded link diagram. We say that $K^\prime$ is the ext. welded link diagram obtained from $K$ by \emph{replacing $K\restriction_{W_K}$} if the following conditions are satisfied:
\begin{enumerate}
\item $K$ and $K^\prime$ are equal in~$N_c$, for all~$c \in C_K$; 
\item $K^\prime$ has no real crossings in $W_{K^\prime}$; 
\item there is a one-to-one correspondence between the arcs of $K\restriction_{W_K}$ and the arcs of~$K^\prime \restriction_{W_{K^\prime}}$, preserving endpoints, orientation and the parity of the numbers of wen marks;
\item there is a one-to-one correspondence between the loops of $K\restriction_{W_K}$ and those of~$K^\prime \restriction_{W_{K^\prime}}$.
\end{enumerate}
Note that an ext. welded link diagram $K^\prime$ has the same Gauss data as $K$ if and only if $K^\prime$ can be deformed through an isotopy of $\Rr^2$ such that it is obtained from $K$ by replacing $K\restriction_{W_K}$.

\begin{lem}
\label{L:Lemma4}
Let $K$ and $K^\prime$ be two ext. welded links diagrams with the same Gauss data. Then $K$ and $K^\prime$ are equivalent. Moreover one can be obtained from the other by isotopy of $\Rr^2$ and a finite sequence of moves of type $(V1)$, $(V2)$, $(V3)$, $(M)$, $(T1)$, and~$(T2)$. 
\end{lem}
\proof
The proof is the same as the one given for the virtual case in~\cite[Lemma~4]{Kamada:Markov}. It should be noted that there is an unsaid difference in the fact that in the extended welded context saying that two diagrams have the same Gauss data implies one condition more than in the virtual and welded cases. This is the condition of preserving wen marks information on the elements of~$G_K$, which is solved by moves $(T1)$ and~$(T2)$.  Compare Definition~\ref{D:Gauss_data} and virtual Gauss data definition in~\cite[Section~4]{Kamada:Markov}. Other proofs for the virtual case can be found in~\cite{Goussarov-Polyak-Viro:2000, Kauffman:Virtual}.
\endproof

We introduce now a particular kind of ext. welded link diagrams. 
\begin{defn}
\label{D:braided_links}
An ext. welded link diagram is said to be \emph{braided} if there exists a point on the plane with respect to which the ext. welded link diagram is braided around. More formally: we can assume the point to be the origin $O$ of $\Rr^2$, and identify $\Rr^2 \setminus \{O\}$ with $\Rr_+ \times S^1$ passing to polar coordinates and considering $S^1$ to be oriented anti-clockwise.  
Let us consider $\pi_2 \colon \Rr_+ \times S^1 \to S^1$ to be the standard projection to the second factor. Then $L$ is a \emph{braided ext. welded link diagram} if:
\begin{enumerate}[label=(\roman*)]
\item $L$ is contained in $\Rr_+ \times S^1$;
\item each component of $L$ is monotone with respect to the coordinate in $S^1$;
\item the restriction $\pi \restriction_{ C_L }$ is injective. 
\end{enumerate}
\end{defn}
A braided ext. welded diagram can be easily represented as the closure of an ext. welded braid diagram, who is uniquely determined up to conjugation. 

For completeness, we quickly recall here the \emph{braiding process} exposed in~\cite{Kamada:Markov} for virtual link diagrams, which can be adapted without modification to ext. welded diagrams. 
Let $L$ be an ext. welded link diagram, and let $N_1, N_2, \ldots, N_n$ be regular neighbourhoods of its real crossings. With an isotopy of $\Rr^2$ we  deform $L$ in such a way so that: all the $N_\ii$s are in $\Rr^2 \setminus \{O\}$; $\pi_2 (N_\ii) \cap \pi_2(N_\jj) = \emptyset$ for $\ii \neq \jj$ in $\{1, \ldots, \nn\}$; each $N_\ii$ contains two oriented arcs, each of which is mapped to $S^1$ by $\pi_2$ homeomorphically with respect to the orientation of~$S^1$. Finally replace $L \restriction{W_L}$ arbitrarily such that the result is a braided ext. welded link diagram, which is equivalent to $L$ (Lemma~\ref{L:Lemma4}).

In the following we will apply the extended Markov moves defined in Subsection~\ref{SS:Extended_Markov_Moves} to braided ext. welded link diagrams in the natural way.  We will still denoted them by~$(M0)$ and~$(M2)$.

\begin{lem}
\label{L:Lemma5}
Let $K$ and $K^\prime$ be braided ext. welded link diagrams, possibly with a different degree as braids, such that $K^\prime$ is obtained from $K$ by replacing $K\restriction_{W_K}$. Then $K$ and $K^\prime$ are related by a finite number of moves of type $(M0)$ and~$(M2)$.
\end{lem}
\proof
The proof is the same as the one given for the virtual case in~\cite[Lemma~5]{Kamada:Markov}, once more with the hidden difference that  arcs or loops of  $K\restriction_{W_K}$ and~$K^\prime\restriction_{W_{K^\prime}}$ could have wen marks. However, the operation of replacing the restrictions for ext. welded diagrams preserves the presence of wen marks, leaving the proof remains unchanged.  
\endproof

\subsection{About the sign reversal move}
\label{SS:About_Sign_Reversal}
Among the extended Markov moves introduced in Subsection~\ref{SS:Extended_Markov_Moves}, the most exotic one is the sign reversal move.  As noted in~\cite[Notation~2.6]{Audoux:2016}, the sign reversal move on a ext. welded link diagram changes its Gauss data by reversing the signs associated to the crossings and leaving everything else unchanged. We refer to the effect of the sign reversal move on the Gauss data as the \emph{sign reversal of the Gauss data}.

\begin{prop}
\label{P:sign_reversal}
The following statements hold:
\begin{enumerate}[label=(\arabic*)]
\item \label{link} For ext. welded link diagrams, sign reversal is a consequence of moves~$(T1) - (T4)$.
\item \label{braided} For braided ext. welded link diagramns, sign reversal is a consequence of moves of type~$(M0)$.
\item For ext. welded braid diagrams, sign reversal is a consequence of moves of type $(M0)$ and $(M1)$. 
\end{enumerate}
\end{prop}
\proof
\begin{enumerate}[label=(\arabic*)]
\item Remark that when move $(T4)$ is applied to a crossing, it locally induces a sign reversal on that crossing. 
Let $L$ be an ext. welded link. It is always possible to introduce a pair of wen marks on each components through a $(T1)$ move. Then, we make one wen mark of each pair slide along the component until it comes back to the starting point, on the other side of the other wen mark of the pair. In such a way it will have passed once on the overstrand of each crossing in which the component is involved ad an overstrand, conjugating the crossings by welded crossings and changing the sign of the crossing.  Then we cancel the pairs of wen marks added through an inverse $(T1)$ move. The resulting diagram $L^\prime$ is the sign reverse of $L$, and is equivalent to $L$ through the moves~$(T1) - (T4)$.
\item Same as point~\ref{link}.
\item This is a direct consequence of point~\ref{braided}, since moves $(M1)$ allow to pass wen marks from the bottom to the top of a diagram. 
\end{enumerate}

\endproof

\begin{lem}
\label{L:Lemma6}
Two braided ext. welded link diagrams with the same Gauss data are related by a finite sequence of Markov moves of types $(M0)$ and~$(M2)$. If the Gauss data of one braided ext. welded link diagram is the sign reverse of the Gauss data of the other, then the two diagrams are related by a finite number of Markov moves of type $(M0)$ and~$(M2)$. 
\end{lem}
\proof
The first part is the same as in \cite[Lemma~6]{Kamada:Markov}. We recall it here for completeness. Let $L$ and $L^\prime$ be braided ext. welded link diagrams with the same Gauss data. Let $c_1, \ldots, c_n$ be the real crossings of~$L$, and $N_1, \ldots, N_n$ the relative regular neighbourhoods. In the same way let $c_1^\prime, \ldots, c_n^\prime$ be the real crossings of $L^\prime$ and $N_1^\prime, \ldots, N_n^\prime$ their neighbourhoods. Let us distinguish two cases, depending on the order of the images of the crossings through $\pi_2$ (Definition~\ref{D:braided_links}).
\begin{enumerate}[label=(\alph*)]
\item \label{It:a} Suppose the crossings $\pi_2(N_1), \ldots, \pi_2(N_n)$ and $\pi_2(N_1^\prime), \ldots, \pi_2(N_n^\prime)$ appear on $S^1$ in the same cyclic order. 
Then, with an isotopy of $\Rr^2$ we can deform $L$ keeping the braidedness, in such a way that $N_i$ and $N_i^\prime$ coincide for $i=1, \ldots, n$, and the restrictions of $L$ and $L^\prime$ to these disks are identical. By Lemma~\ref{L:Lemma5}, $L$ and $L^\prime$ with a finite number of Markov moves of types $(M0)$ and~$(M2)$.
\item \label{it:b} Suppose that $\pi_2(N_1), \ldots, \pi_2(N_n)$ and $\pi_2(N_1^\prime), \ldots, \pi_2(N_n^\prime)$ do not appear on $S^1$ in the same cyclic order. It is enough to treat the case when only a pair, for example $\pi_2(N_1)$ and $\pi_2(N_2)$ are exchanged. Then, with a finite sequence of moves~$(M_0)$ and $(M_1)$ (details in \cite[Lemma~6]{Kamada:Markov}) it is possible to move one crossing, preserving at each step the braidedness and the Gauss data, in order to reconduct ourselves to case~\ref{It:a}.
\end{enumerate}

Let us now consider the case in which $L$ and $L^\prime$ are braided ext. welded link diagrams with sign reversed Gauss data. We consider two cases as before. 
\begin{enumerate}[label=(\alph*)]
\item \label{It:2a} Suppose the crossings $\pi_2(N_1), \ldots, \pi_2(N_n)$ and $\pi_2(N_1^\prime), \ldots, \pi_2(N_n^\prime)$ appear on $S^1$ in the same cyclic order. Again with an isotopy of $\Rr^2$ we can deform $L$ keeping the braidedness, in such a way that $N_i$ and $N_i^\prime$ coincide for $i=1, \ldots, n$, but the restrictions of $L$ and $L^\prime$ to these disks present opposite crossings. Let us apply to $L$ a sign reversal. Then we obtain a braided ext. welded diagram $L^{\prime \prime}$ which is equivalent to $L$ via $M0$ moves by Proposition~\ref{P:sign_reversal}, whose crossings regular neighbourhoods $N_1^{\prime\prime}, \ldots, N_n^{\prime\prime}$ contain a real crossing conjugated by welded crossings. With an isotopy of $\Rr^2$, push the welded crossings outside of the regular neighbourhood, deforming  $L^{\prime\prime}$ in such a way that $N^{\prime\prime}_i$ and $N_i^\prime$ coincide for $i=1, \ldots, n$, and the restrictions of $L^{\prime\prime}$ and $L^\prime$ to these disks are identical. Then by Lemma~\ref{L:Lemma5}, one can pass from $L^{\prime\prime}$ to $L^\prime$ with a finite number of Markov moves of types $(M0)$ and~$(M2)$.
\item Suppose the crossings $\pi_2(N_1), \ldots, \pi_2(N_n)$ and $\pi_2(N_1^\prime), \ldots, \pi_2(N_n^\prime)$ appear on $S^1$ in the same cyclic order. Then with the manoeuvre recalled in point~\ref{it:b} of the first part, we reconduct this case to point~\ref{It:2a} of this part. 
\end{enumerate} 
\endproof

The following is a corollary of the first part of Lemma~\ref{L:Lemma6}, which is a direct consequence of the fact that the braiding process does not change the Gauss data. In the case of virtual link diagrams, it appears in~\cite[Corollary~7]{Kamada:Markov}.

\begin{cor}
\label{C:Corollary7}
For an ext. welded link diagram $K$, a braided ext. welded link diagram obtained by the braiding process is uniquely determined up to ext. welded Markov moves $(M0)$ and~$(M2)$. 
\end{cor}

\section{A Markov theorem for welded extended diagrams}
\label{S:MarkovTheorem}
\begin{thm}
\label{T:Markov}
Two ext. welded braid diagrams that admit closure have equivalent closures as ext. welded link diagrams if and only if they are related by a finite sequence on the following moves:
\begin{description}[before={\renewcommand\makelabel[1]{\bfseries ##1}}]
\item[$(M0)$] \label{itm:M0} isotopy of $\Rr^2$ and generalized  Reidemeister moves;
\item[$(M1)$] \label{itm:M1} conjugation in the welded ext. braid group~$\WBE\nn$;
\item[$(M2)$] \label{itm:M2} a right stabilization of positive, negative or welded type, and its inverse operation.
\end{description}
\end{thm}
\proof
Let $b$ and $b^\prime$ be ext. welded braid diagrams that admit closure related by $(M0)$ moves. Then their closures $\widehat{b}$ and $\hat{b^\prime}$ are equivalent ext. welded link diagrams for definition of generalized Reidemeister equivalence.
Suppose that  $b$ and $b^\prime$ are related by moves of type $(M1)$ and~$(M2)$: then also in these cases their closures $\hat{b}$ and $\hat{b^\prime}$ are clearly equivalent as ext. welded link diagrams. 

\begin{figure}[htb]
\centering
\includegraphics[scale=0.65]{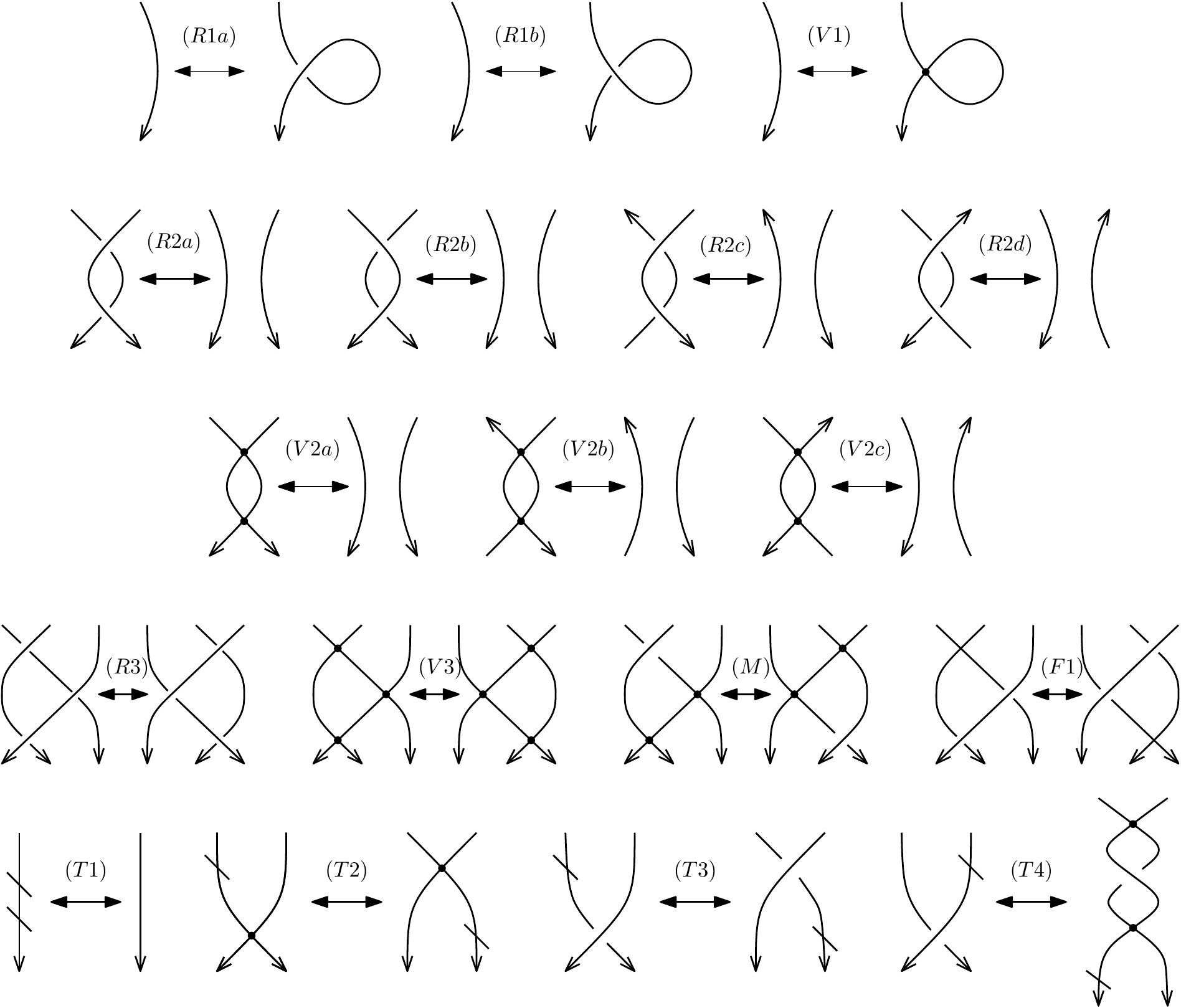}
\caption{Oriented generalized Reidemeister moves.}
\label{F:Oriented_Reidemeister}
\end{figure}

On the other hand, let $K$ and $K^\prime$ be ext. welded link diagrams representing the same exteded welded link. Then there is a finite sequence of ext. welded link diagrams $K=K_0, K_1, \ldots, K_s=K^\prime$ such that $K_i$ is obtained from $K_{i-1}$ by oriented generalized Reidemeister moves of kind $(R1a)$, $(R1b)$, $(V1)$, $(R2a)$, $(R2b)$, $(R2c)$, $(R2d)$, $(V2a)$, $(V2b)$, $(V2c)$, $(R3)$, $(V3)$, $(M)$, $(F1)$, $(T1)$, $(T2)$, $(T3)$, or $(T4)$, as shown in Figure~\ref{F:Oriented_Reidemeister}. 
Applying the braiding process to each $K_i$, we obtain a braided ext. welded link $\widetilde{K}_i$ with the same Gauss data as $K_i$.
By Corollary~\ref{C:Corollary7}, $\widetilde{K}_i$ is uniquely determined up to moves $(M0)$ and $(M2)$. 
To prove the statement it is enough to check that for each $i=1, \ldots, s$, $\widetilde{K}_i$ and $\widetilde{K}_{i-1}$ are equivalent up to moves $(M0) - (M2)$. 
All the moves except $(T1)$, $(T2)$, $(T3)$ and $(T4)$ are considered in~\cite[Theorem~2]{Kamada:Markov}.
So let us consider the remaining cases. 
Suppose that $K_i$ is obtained by $K_{i-1}$ by a $(T1)$, $(T2)$, $(T3)$ or $(T4)$ move. Then, let $\Delta$ be a $2$-disk in $\Rr^2$ that contains one of these moves, and let $\Delta^c$ be such that $K_i \cap \Delta^c = K_{i-1} \cap \Delta^c$. 
Deform $K_i$ and $K_{i-1}$ by an isotopy of $\Rr^2$ in such a way that $K_i \cap \Delta$ and $K_{i-1} \cap \Delta$ satisfy the condition to be a braided ext. welded link diagram.
Applying the braiding process to  $K_i \cap \Delta^c$ and $K_{i-1} \cap \Delta^c$ we obtain diagramd $\widetilde{K}_i^\prime$ and $\widetilde{K}_{i-1}^\prime$ such that:
\begin{equation*}
\widetilde{K}^\prime_i \cap \Delta= {K_i} \cap \Delta  \quad \mbox{and} \quad \widetilde{K}^\prime_{i-1} \cap \Delta = K_{i-1} \cap \Delta  \quad \mbox{and} \quad \widetilde{K}^\prime_i \cap \Delta^c = \widetilde{K}^\prime_{i-1} \cap \Delta ^c
\end{equation*}
Then $\widetilde{K}_i^\prime$ and $\widetilde{K}_{i-1}^\prime$ are related by a $(M0)$ move corresponding to $(T1)$, $(T2)$, $(T3)$, or~$(T4)$. 
Since  $\widetilde{K}_i^\prime$  has the same Gauss data as $K_i$, it is equivalent by Markov moves to  $\widetilde{K}_i$ (Lemma~\ref{L:Lemma6}). Same for  $\widetilde{K}_{i-1}^\prime$ and~$\widetilde{K}_i$. Therefore $\widetilde{K}_i$ and $\widetilde{K}_{i-1}$ are equivalent by Markov moves.
\endproof

\section{Extended welded knots and horizontal mirror images}
\label{S:Mirror}

\begin{figure}[htb]
\centering
\includegraphics[scale=0.6]{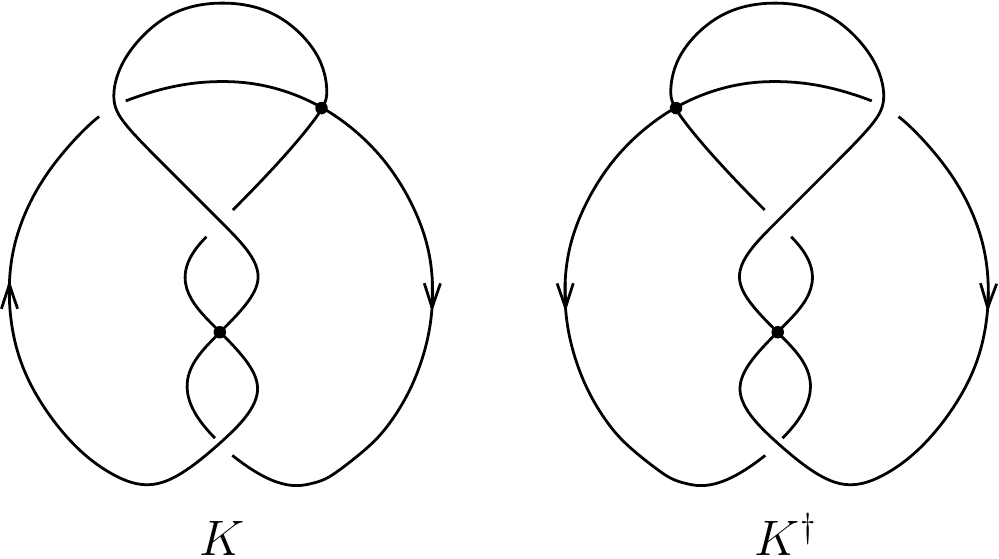}
\caption{An extened welded knot diagram $K$ and its horizontal mirror image $K^\dagger$.}
\label{F:MirrorImage}
\end{figure}

Let $K$ be an ext. welded diagram. Its \emph{horizontal mirror image} $K^\dagger$ is its reflection with respect to a line on the plane of the diagram, as in Figure~\ref{F:MirrorImage}. In this section we show that ext. welded knots are equivalent to their horizontal mirror image. In particular, we show that the horizontal mirror image $K^\dagger$ of an ext. knot diagram~$K$ is equivalent to the sign reversal of~$K$.

\begin{prop}
\label{P:Mirror}
Every ext. welded knot diagram is equivalent to its horizontal mirror image.
\end{prop}

\proof
Consider an ext. welded knot diagram~$K$. Let us denote by~$sK$ the sign reversal of $K$. It is easy to see that $K^\dagger$ has the same Gauss data as~$sK$ (also shown in~\cite{Ichimori-Kanenobu:2012}). By Lemma~\ref{L:Lemma4}, $sK$ and $K^\dagger$ are equivalent. Then, $K^\dagger$ is equivalent to $K$ by Proposition~\ref{P:sign_reversal}. 
\endproof

\section*{Acknowledgements}
The author thanks Seiichi Kamada for helpful comments. This work is part of a joint research project with Seiichi Kamada, supported by JSPS KAKENHI Grant Number 16F16793.
During the writing of this paper, the author was supported by a JSPS Postdoctral Fellowship For Foreign Researchers.

\bibliographystyle{alpha}
\bibliography{Markov_Diagram.bib}{}

\end{document}

%% file: Macro.tex
\newcommand\Zz{\mathbb{Z}}

\newcommand\Rr{\mathbb{R}}

\newcommand\WB[1]{WB_{#1}}	
\newcommand\WBE[1]{WB_{#1}^{ext}}	
\newcommand\LB[1]{LB_{#1}}	
\newcommand\rB[1]{rB_{#1}}	

\newcommand\ii{i}

\newcommand\jj{j}

\newcommand\nn{n}
\newcommand\nno{{n-1}}

\makeatletter
\newcommand{\oset}[2]{%
  {\mathop{#2}\limits^{\vbox to -.5\ex@{\kern-\tw@\ex@
   \hbox{\scriptsize #1}\vss}}}}
\makeatother

\newcommand\rr[1]{\rho_{\hspace{-0.3ex}#1}^{\null}}	
\newcommand\rrq[2]{\rho_{\hspace{-0.3ex}#1}^{#2}}	

\newcommand\sig[1]{\sigma_{\hspace{-0.3ex}#1}^{\null}} 
\newcommand\sigg[2]{\sigma_{\hspace{-0.3ex}#1}^{#2}}	
\newcommand\siginv[1]{\sigg{#1}{-1}}	
